\documentclass[12pt,leqno,a4paper]{article}
\usepackage{amsfonts}
\usepackage{amsmath, amssymb, amsthm, amsfonts}

\numberwithin{equation}{section}

\begin{document}
\title{On the error term in Weyl's law for the Heisenberg manifolds(II)}
\author{  ZHAI Wenguang \\ \\
Science in China series A: Mathematics, Vol. 52(2009),857-874}
\date{}

\footnotetext[0]{2000 Mathematics Subject Classification: 11N37,
35P20, 58J50.} \footnotetext[0]{Key Words: Heisenberg manifold,
Weyl's law, error term, mean square.} \footnotetext[0]{This work is
supported by  National Natural Science Foundation of China(Grant No.
10771127).} \maketitle

{\bf Abstract.} In this paper we study the mean square of  the error
term in the Weyl's law of  an irrational $(2l+1)$-dimensional
Heisenberg manifold . An asymptotic formula is established.

{\section{Introduction}}

Let $(M, g)$ be a closed $n$-dimensional Riemannian manifold with
metric $g$ and Laplace-Beltrami operator $\Delta.$ Let $N(t)$ denote
its spectral counting function , which is defined as the number of
the eigenvalues of $\Delta$ not exceeding $t.$ H\"ormander \cite{Ho}
proved that the Weyl's law
\begin{equation}
N(t)=\frac{vol(B_n)vol(M)}{(2\pi)^n}t^{n/2}+O(t^{(n-1)/2})
\end{equation}
holds, where $vol(B_n)$ is the volume of the $n$-dimensional unit
ball.

Let $$R(t): =N(t)-\frac{vol(B_n)vol(M)}{(2\pi)^n}t^{n/2}.$$
H\"ormander's estimate (1.1) in general is sharp , as the well-known
example of the sphere $S^n$ with its canonical metric shows
\cite{Ho}. However, it is a very difficult problem to determine the
optimal bound of $R(t)$ in any given manifold, which depends on the
properties of the associated geodesic flow. Many improvements have
been obtained for certain types of manifolds, see \cite{BG, BE,
BL,CPT, FR, G, Hu, IVr, KP, PT, Vo}.

\subsection{\bf The Weyl's law for ${\Bbb T}^2$: the Gauss circle problem}\

The simplest compact manifold with integrable geodesic flow is the
$2$-torus ${\Bbb T}^2={\Bbb R}^2/{\Bbb Z}^2.$ The exponential
functions $e(mx+ny)(m,n\in{\Bbb Z})$ form a basis of eigenfunctions
of the Laplace operator $\Delta=\partial_x^2+\partial_y^2,$ which
acts on functions on ${\Bbb T}^2.$ The corresponding eigenvalues are
$4\pi^2(m^2+n^2), m,n\in {\Bbb Z}.$ The spectral counting function
$$N_I(t)=\{\lambda_j\in Spec(\Delta): \lambda_j\leq t\}$$
is equal to the number of lattice points of ${\Bbb Z}^2$ inside a
circle of radius $\sqrt t/2\pi.$ The well-known Gauss circle problem
is to study the properties of the error term of the function
$N_I(t).$

In this case , the formula (1.1) becomes
\begin{equation}
N_I(t)=\frac{t}{4\pi}+O(t^{1/2}),
\end{equation}
which is the classical result of Gauss.  Let $R_I(t)$ denote the
error term in (1.2). Many authors improved the upper bound estimate
of $R_I(t).$ The latest result is due to Huxley\cite{Hu}, which
reads
\begin{eqnarray}
R_I(t)\ll t^{131/416}\log^{26947/8320} t.
\end{eqnarray}
Hardy\cite{Ha} conjectured that
\begin{equation}
R_I(t)\ll t^{1/4+\varepsilon},
\end{equation}
which is supported by the asymptotic formula
\begin{eqnarray*}
 \int_1^T|R_I(t)|^2dt= \frac{1}{6\pi^3}\sum_{n=1}^\infty\frac{r^2(n)}{n^{3/2}}T^{3/2}+O(T\log^{2} T)
\end{eqnarray*}
proved in \cite{Ka} , where $r(n)$ denotes the number of ways $n$
can be written as a sum of two squares.

Tsang\cite{Ts}  first   proved that    the asymptotic formula
\begin{equation}
\int_1^T R_I^k(t)dt= c_kT^{1+k/4} +O(T^{1+k/4-\delta_k+\varepsilon})
\end{equation}
holds for $k=3$ and $k=4$ with $\delta_3=1/14$ and $\delta_4=1/23,$
where $c_k(k\geq 3)$ and $\delta_k>0(k\geq 3)$ are explicit
constants.

In \cite{Z1}, the author proved that    (1.5) holds for any integer
$3\leq k\leq 9$  . When $k=4$, in \cite{Z2} the author proved that
we can take $\delta_4=3/28$ in (1.5).

\subsection{\bf The Weyl's law for $(2l+1)$-dimensional rational  Heisenberg
manifold}\

Let $l\geq 1$ be a fixed integer and  $(H_l/\Gamma, g)$ be a
$(2l+1)$-dimensional Heisenberg manifold with a metric $g$. When
$l=1$, in \cite{PT} Petridis and Toth proved that
$R(t)=O(t^{5/6}\log t)$ for a special metric. Later in \cite{CPT}
this bound was improved to $O(t^{119/146+\varepsilon})$ for all
left-invariant Heisenberg metrics. For $l>1$   Khosravi and
Petridis\cite{KP} proved that $ R(t)=O(t^{l-7/41})$ holds for
rational Heisenberg manifolds. Both in \cite{CPT} and \cite{KP},
they first established a $\psi$-expression of $R(t)$ and then used
the van der Corput method of exponential sums  .
 If substituting Huxley's result of \cite{Hu} into the arguments of \cite{CPT} and \cite{KP}, we
 can get that the estimate
\begin{equation}
R(t)=O(t^{l-77/416}(\log t)^{26947/8320})
\end{equation}
holds for all rational $(2l+1)$-dimensional Heisenberg manifolds,
which corresponds to Huxley's result (1.3).

It was conjectured that for rational Heisenberg manifolds,  the
pointwise estimate
\begin{equation}
R(t)\ll t^{l-1/4+\varepsilon}
\end{equation}
holds, which was proposed  in  Petridis and Toth \cite{PT} for the
case $l=1$ and in Khosravi and Petridis\cite{KP} for the case $l>1.$
As an evidence of this conjecture, Petridis and Toth proved the
following $L^2$ result for $H_1$
\begin{eqnarray*}
\int_{I^3}\left|N(t;u)-\frac{1}{6\pi^2}vol(M(u))t^{3/2}\right|^2du\leq
C_\delta t^{3/2+\delta},
\end{eqnarray*}
where $I=[1-\varepsilon,1+\varepsilon].$ They also proved
$$\frac{1}{T}\int_T^{2T}\left|N(t)-\frac{1}{6\pi^2}vol(M)t^{3/2}\right|dt\gg T^{3/4}.$$

Let $M=(H_l/\Gamma, g_l)$ be a $(2l+1)$-dimensional Heisenberg
manifold with the metric
 \[g_l:=\left(
\begin{array}{llcl}
 I_{2l\times 2l}&0\\
0&2\pi\\
\end{array}
\right),\] where $I_{2l\times 2l}$ is the identity matrix.

M. Khosravi and John A. Toth\cite{KT} proved that
\begin{equation}
\int_1^TR^2(t)dt=C_{2,l}T^{2l+1/2}+O(T^{2l+1/4+\varepsilon}),
\end{equation}
where $C_{2,l}$ is an explicit constant .

 M. Khosravi \cite{K} proved that  the asymptotic formula
\begin{equation}
\int_1^TR^3(t)dt=C_{3,l}T^{3l+1/4}+O(T^{3l+3/14+\varepsilon})
\end{equation}
is true for some explicit constant $C_{3,l}.$

Recently, the author\cite{Z4} proved that the asymptotic formula
\begin{equation}
\int_1^T R^k(t)dt= C_{k,l}T^{k(l-1/4)+1}
+O(T^{k(l-1/4)+1-\eta_k+\varepsilon})
\end{equation}
holds for any $3\leq k\leq 9,$ where $C_{k,l}$ and $\eta_k>0$ are
explicit constants. Especially (1.10) holds for $k=3$ with
$\eta_3=1/4$ and for $k=4$ with $\eta_4=3/28.$
\bigskip

The moments  problem of $R(t)$ becomes very difficult for the
irrational Heisenberg manifolds  even  when we study only the mean
square case. For the definition of rationality of Heisenberg
manifolds, see\cite{KP}. The aim of this paper is to study the mean
square of the error term in the Weyl's law for the
$(2l+1)$-dimensional irrational Heisenbergs.

 The plan of this paper is as follows.  In Section 2 we shall state
our main results. In Section 3 we state some background of the
Heisenberg manifolds and give a  $ \psi$-expression of $R(t).$ In
Section 4 are some preliminary Lemmas. We shall prove our theorem in
Section 5.

{\bf Notations.} For a real number  $t,$ let $[t]$ denote the
integer part of $t,$ $\{t\}=t-[t],$ $\psi(t)=\{t\}-1/2,$  $\Vert
t\Vert=\min(\{t\},1-\{t\}),$ $e(t)=e^{2\pi it}.$ $\varepsilon$
always denotes a sufficiently small positive constant. ${\Bbb C},
{\Bbb R},
 {\Bbb Z},  {\Bbb N}$ denote the set of complex numbers, the set of real numbers, the set of
integers, the set of positive integers, respectively. $n\sim M$
means that $N<n\leq 2N$ and $n\asymp N$ means $c_1N\leq n\leq c_2N$
for some positive constants $0<c_1<c_2.$     $SC(\Sigma)$ denotes
the summation condition of the sum $\Sigma.$ Throughout this paper ,
${\cal L}$ always denotes $\log T.$

\section{\bf Main results}

From now on, we always suppose that $R(t)$ denote the error term in
the Weyl's law for the $(2l+1)$-dimensional Heisenberg manifold with
the irrational metric
\[g_l(\theta):=\left(
\begin{array}{llcl}
 I_{2l\times 2l}&0\\
0& 2\pi/\theta\\
\end{array}
\right),\] where $\theta>0$ is an irrational number.

The mean square of $R(t)$ is closely related to the {\it
approximation type} of $\theta.$ We  recall a few facts from the
theory of Diophantine approximation: by the  {\it approximation
type} $\gamma(\alpha)$ of an irrational real number $\alpha$ we
denote the infimum of all reals $r$ for which there exists a
constant $c(r,\alpha, \epsilon)$ such that the inequality
\begin{equation}
|\alpha-p/q|\geq c(r,\alpha, \varepsilon)q^{-r-1-\varepsilon},\
 \end{equation}
for any $p\in {\Bbb Z}$ and any $q\in {\Bbb N}.$ Obviously
$\gamma(\alpha)\geq 1$ for all irrationals. By Roth's theorem
\cite{R}, if $\alpha$ is an algebraic irrational, then
$\gamma(\alpha)=1.$ Furthermore, $\gamma(\alpha)=1$ for almost all
irrationals(see Khinchin\cite{KH}).

 {\bf Theorem 1.} Suppose $\theta>0$ is an irrational number
of finite {\it type} $\gamma .$ Then we have
\begin{eqnarray}
\int_1^TR^2(x)dx=\frac{2^{9/2-4l}{\cal
C}_{l,\theta}}{(4l+1)(l-1)!^2\pi^{2l+3/2}}T^{2l+1/2} +O(T^{2l
+\frac{4\gamma+1}{8\gamma+4}+\varepsilon}),
\end{eqnarray}
where
$${\cal C}_{l,\theta}: =\sum_{h=1}^\infty\sum_{r>\theta h }\frac{(1-\frac{\theta
h}{2r-\theta h})^{2l-2}}{h^{1/2}(2r-\theta h)^{5/2}}.$$

{\bf Corollary.} If $\theta$ is an algebraic irrational, then we
have
\begin{eqnarray}
\int_1^TR^2(x)dx=\frac{2^{9/2-4l}{\cal
C}_{l,\theta}}{(4l+1)(l-1)!^2\pi^{2l+3/2}}T^{2l+1/2} +O(T^{2l
+\frac{5}{12}+\varepsilon}).
\end{eqnarray}
Furthermore  (2.3) holds for almost all  irrational  number
$\theta>0.$

For the $3$-dimensional case, we have the following more general
Theorem 2. We omit its proof  since it is almost    the same  as
that of Theorem 1.

{\bf Theorem 2. } Let $R(t)$ denote the error term in the Weyl's law
for the $3$-dimensional Heisenberg manifold with the metric
 \[g:=\left(
\begin{array}{llcl}
 h&&0\\
 &0&g_3\\
\end{array}
\right),\] where $g_3>0, $ $h=(h_{ij})\ (i,j=1,2)$ with
$h_{11}h_{22}-h_{12}^2>0$.   Suppose that $\pm id^2(d>0)$ are the
eigenvalues of the matrix $h^{-1}J,$ where $J$ is the standard
symplectic $2\times 2$ matrix. Let $\theta=2\pi/g_3d^2.$

If $\theta$  is an irrational number of finite {\it approximation
type} $\gamma,$ then we have
\begin{eqnarray}
\int_1^TR^2(x)dx=\frac{2^{1/2}{\cal C}_{1,\theta}}{
5d^3\pi^{7/2}}T^{5/2} +O(T^{
\frac{20\gamma+9}{8\gamma+4}+\varepsilon}).
\end{eqnarray}

{\bf Remark.} If the value $\theta$ in Theorem 2 is a rational
number, then the error term $O(T^{
\frac{20\gamma+9}{8\gamma+4}+\varepsilon})$ in (2.4) can be replaced
by $O(T^{9/4+\varepsilon}).$

\section{\bf Background of Heisenberg manifolds and the $\psi$-expression of $R(t)$}

In this section, we first review  some background of the Heisenberg
manifolds. The reader can see \cite{Fo}, \cite{GW} , \cite{St} for
more details. Finally, we give an $\psi$-expression of $R(t)$.

\subsection{\bf Heisenberg manifolds}\

 Suppose  $x\in {\Bbb R}^l$ is a row vector  and $y\in {\Bbb R}^l$ is a column vector. Define
 \[\gamma(x,y,t)=\left(
\begin{array}{llcl}
 1&x&t\\
0&I_l&y\\
0&0&1\\
\end{array}
\right), \ \ X(x,y,t)=\left(
\begin{array}{llcl}
 0&x&t\\
0&0&y\\
0&0&0\\
\end{array}
\right).\] The $(2l+1)$-dimensional Heisenberg group $H_l$ is
defined by
$$H_l=\{\gamma(x,y,t):  x, y\in {\Bbb R}^l, t\in {\Bbb R}\},$$
its Lie algebra is
$$\mathfrak{H}_l=\{X(x,y,t):  x, y\in {\Bbb R}^l, t\in {\Bbb R}\}.$$
We say $\Gamma$ is uniform discrete subgroup of $H_l$ if
$H_l/\Gamma$ is compact. A $(2l+1)$-dimensional Heisenberg manifold
is a pair $(H_l/\Gamma, g)$ for which $\Gamma$ is a uniform discrete
subgroup of $H_l$   and $g$ is a left $H_l$-invariant metric.

For every  $r$-tuple $(r_1,r_2,\cdots,r_l)\in {\Bbb N}^l$ such that
$r_j|r_{j+1}\ (j=1,2,\cdots,l-1)$, let $r{\Bbb Z}^l$ denote the
$l$-tuples $x=(x_1,x_2,\cdots,x_l)$ with $x_j\in r_j{\Bbb Z}$.
Define
$$\Gamma_r=\{\gamma(x,y,t): x\in r{\Bbb Z}^l, y\in r{\Bbb Z}^l, t\in {\Bbb Z}\}.$$
It is clear that $\Gamma_r$ is a  uniform discrete subgroup of
$H_l$.  According to Theorem 2.4 of \cite{GW}, the subgroup
$\Gamma_r$ classifies all the uniform discrete subgroups of $H_l$ up
to automorphisms. Thus (see \cite{GW}, Corollary 2.5) given any
Riemannian Heisenberg manifold $M=(H_l/\Gamma, g)$, there exists a
unique $l$-tuple $r$ as before and a left-invariant metric ${\tilde
g}$ on $H_l$ such that $M$ is isometric to $(H_l/\Gamma, {\tilde
g}).$ So (see \cite{GW}, 2.6(b)) we can replace the metric $g$ by
$\phi^*g,$ where $\phi$ is an inner automorphism such that the
direct sum split of the Lie algebra $\mathfrak{H}_l={\Bbb
R}^{2l}\oplus\mathfrak{Z}$ is orthogonal . Here  $\mathfrak{Z}$ is
the center of the Lie algebra and
 $${\Bbb R}^{2l}=\left\{\left(
\begin{array}{llcl}
 0&x&0\\
0&0&y\\
0&0&0\\
\end{array}
\right): x,y\in {\Bbb R}^l\right\}.$$ With respect to this
orthogonal split of $H_l$ the metric $g$ has the form \[\left(
\begin{array}{llcl}
 h& 0\\
0&g_{2l+1}\\
 \end{array}
\right),\] where $h$ is a positive-definite $2l\times 2l$ matrix and
$g_{2l+1}>0$ is a real number.

The volume of the Heisenberg manifold is given by
$$vol(H_l/\Gamma,g)=|\Gamma_r|\sqrt{det(g)}$$
with $|\Gamma_r|=r_1r_2\cdots r_l$ for $r=(r_1,r_2,\cdots, r_l).$

\subsection{\bf The spectrum of Heisenberg manifolds and  the $\psi$-expression of $R(t)$}\

Let $\Sigma$ be the spectrum of the Laplacian on $M=(H_l/\Gamma,
g_l(\theta)),$ where the eigenvalues are counted with
multiplicities. According to \cite{GW}(P. 258), $\Sigma$ can be
divided into two parts $\Sigma_I$ and $\Sigma_{II},$ where
$\Sigma_I$ is the spectrum of $2l$-dimensional torus
 and $\Sigma_{II}$ contains all eigenvalues of the form
 $$\frac{4\pi^2m^2}{g_{2l+1}}+\sum_{j=1}^{l}2\pi m(2n_j+1),\ \ m\in {\Bbb N}, n_j\in {\Bbb N}\cup \{0\},$$
 each eigenvalue counted with the multiplicity  $2m^l.$

We have   the following $\psi$-expression of $R(t).$

{\bf Lemma 3.1.} We have
\begin{eqnarray}
&&R(2\pi x)=-\frac{4}{2^l(l-1)!}\sum_{1\leq m \leq \sqrt
\frac{x}{\theta}}m(x-\theta m^2
)^{l-1}\psi\left(\frac{x}{2m}-\frac{\theta m}{2}-\frac{l}{2}\right)
\end{eqnarray}
$$+O(x^{l-1/2}).$$

In \cite{Z4}, the author proved Lemma 3.1 when $\theta=1.$ However,
the proof for the general case is almost the same. So we omit the
details of the proof.

\section{Some preliminary Lemmas}

We need the following Lemmas. Lemma 4.1 is due to Vaaler\cite{Va}.
Lemma 4.2  is well-known; see for example, Heath-Brown\cite{He2}.
Lemma 4.3 is Theorem 2.2 of Min\cite{Mi},  see also Lemma 6 of
Chapter 1 in \cite{V}. A weaker version of Lemma 4.3 can be found in
\cite{Ku}, which also suffices for our proof. Lemma 4.4 and Lemma
4.6 provide several estimates
   about the quantity
 $$\alpha(\theta;h_1,h_2,n_1,n_2):=\sqrt{h_1(2n_1-\theta h_1)}-\sqrt{h_2(2n_2-\theta h_2)},$$
which play an essential role in our proof.

{\bf Lemma 4.1.} Let $H\geq 2$ be any real number. Then
$$\psi(u)=\sum_{1\leq |h|\leq H}a(h)e(hu)+O(\sum_{0\leq |h|\leq H}b(h)e(hu)),$$
where $a(h)$ and $b(h)$ are functions such that $a(h)\ll 1/|h|,
b(h)\ll 1/H.$

{\bf Lemma 4.2.}   Let $H\geq 2$ be any real number. Then
$$\psi(u)=-\sum_{1\leq |h|\leq H}\frac{e(hu)}{2\pi ih}+O\left(\min(1,\frac{1}{H\Vert u\Vert}) \right).$$

{\bf Lemma 4.3.} Suppose  $A_1,\cdots, A_5$ are absolute positive
constants, $f(x)$ and $ g(x)$ are algebraic functions in $[a,b]$ and
\begin{eqnarray*}
&&\frac{A_1}{R}\leq |f^{''}(x)|\leq\frac{A_2}{R},\  \ \ |f^{'''}(x)|\leq\frac{A_3}{RU},\ \ U\geq 1, \\
&&|g(x)|\leq A_4G,\ \ \  |g^{'}(x)|\leq A_5GU_1^{-1},\ \ U_1\geq 1,
\end{eqnarray*}
$[\alpha,\beta]$ is the image of  $[a,b]$ under the mapping
$y=f^{'}(x)$, then
\begin{eqnarray*}
\sum_{a<n\leq b}g(n)e(f(n))&=&e^{\pi i/4}\sum_{\alpha<u\leq
\beta}b_u
\frac{g(n_u)}{\sqrt{f^{''}(n_u)}}e\left(f(n_u)-un_u \right)\\
&&+O\left(G\log(\beta-\alpha+2)+G(b-a+R)(U^{-1}+U_1^{-1})\right)\\
&&+ O\left(G\min\left[\sqrt R,
\max\left(\frac{1}{<\alpha>},\frac{1}{<\beta>}\right)\right]\right),
\end{eqnarray*}
where  $n_u$  is the solution of $f^{'}(n)=u$,
\begin{eqnarray*}
<t>=\left\{\begin{array}{ll}
\Vert t\Vert,&\mbox{if $t$ not an integer,}\\
\beta-\alpha,& \mbox{if $t$ an integer,}
\end{array}\right.
\end{eqnarray*}
\begin{eqnarray*}
b_u=\left\{\begin{array}{ll}
1 ,&\mbox{if $\alpha<u<\beta$, or $\alpha, \beta$ not integers ,}\\
1/2,& \mbox{if $\alpha$ or $\beta$ are integers,}\\
\end{array}\right.
\end{eqnarray*}
\begin{eqnarray*}
\sqrt{f^{\prime\prime}}=\left\{\begin{array}{ll}
 \sqrt{f^{\prime\prime}},&\mbox{if $ f^{\prime\prime}>0,$}\\
 i\sqrt{|f^{\prime\prime}|},& \mbox{if $ f^{\prime\prime}<0. $}
\end{array}\right.
\end{eqnarray*}

{\bf Lemma 4.4.} Suppose $\theta>0$ is an irrational number,
$H_1\geq 2, H_2\geq 2, N_1\geq 2, N_2\geq 2, \Delta>0.$ Let ${\cal
A}_{\theta}(H_1, H_2, N_1, N_2;\Delta)$ denote the number of
solutions of the inequality
\begin{equation}
|\sqrt{h_1(2n_1-\theta h_1)}-\sqrt{h_2(2n_2-\theta h_2)}|\leq \Delta
\end{equation}
for which
$$h_1\sim H_1, h_2\sim H_2, n_1\sim N_1, n_2\sim N_2, n_1>\theta h_1, n_2>\theta h_2,$$
then
\begin{eqnarray}
{\cal A}_{\theta}(H_1, H_2, N_1, N_2;\Delta)&&\ll \Delta
(H_1H_2N_1N_2)^{3/4}\\
&&\ \ \ \ +(H_1H_2N_1N_2)^{1/2}\log^2 H_1H_2N_1N_2.\nonumber
\end{eqnarray}

\begin{proof}
If $\Delta\geq (H_1H_2N_1N_2)^{1/4}/100,$ then trivially we have
$${\cal A}_{\theta}(H_1, H_2, N_1, N_2;\Delta)\ll H_1H_2N_1N_2\ll
\Delta (H_1H_2N_1N_2)^{3/4}.$$

Now suppose $\Delta< (H_1H_2N_1N_2)^{1/4}/100.$ In this case $
H_1N_1\asymp H_2N_2.$ Without loss of generality, suppose $H_1\leq
H_2,$ then $N_1\gg N_2.$ It is easy to see that if (4.1) holds ,
then
\begin{eqnarray*}
|2(h_1n_1-h_2n_2)-\theta(h_1^2-h_2^2)|&&\leq
\Delta(\sqrt{h_1n_1}+\sqrt{h_2n_2})\\
&&\leq 2\Delta(\sqrt{H_1N_1}+\sqrt{H_2N_2})\\
&&\leq C\Delta \sqrt[4]{H_1N_1H_2N_2}
\end{eqnarray*}
for some $C>0.$ Thus we have
$$\frac{h_2n_2}{h_1}+\frac{\theta(h_1^2-h_2^2)}{2h_1}-\frac{C\Delta  (H_1N_1H_2N_2)^{\frac{1}{4}}}{2h_1}<n_1\leq
 \frac{h_2n_2}{h_1}+\frac{\theta(h_1^2-h_2^2)}{2h_1}+\frac{C\Delta  (H_1N_1H_2N_2)^{\frac{1}{4}}}{2h_1},$$
which implies that
\begin{eqnarray}
&&{\cal A}_{\theta}(H_1, H_2, N_1, N_2;\Delta)\\
&&\ll \sum_{h_1\sim H_1}\sum_{h_2\sim H_2}\sum_{n_2\sim
N_2}\left[\frac{h_2n_2}{h_1}+\frac{\theta(h_1^2-h_2^2)}{2h_1}+\frac{C\Delta
(H_1N_1H_2N_2)^{\frac{1}{4}}}{2h_1}\right]\nonumber\\&& \ \ \ \ \ \
\ \ \ \ -\sum_{h_1\sim H_1}\sum_{h_2\sim H_2}\sum_{n_2\sim
N_2}\left[\frac{h_2n_2}{h_1}+\frac{\theta(h_1^2-h_2^2)}{2h_1}-\frac{C\Delta
(H_1N_1H_2N_2)^{\frac{1}{4}}}{2h_1}\right]\nonumber\\
&&=\Sigma_1+\Sigma_2-\Sigma_3,\nonumber
\end{eqnarray}
where
\begin{eqnarray*}
&&\Sigma_1=\sum_{h_1\sim H_1}\sum_{h_2\sim H_2}\sum_{n_2\sim
N_2}\frac{C\Delta (H_1N_1H_2N_2)^{\frac{1}{4}}}{h_1},\\
&&\Sigma_2=\sum_{h_1\sim H_1}\sum_{h_2\sim H_2}\sum_{n_2\sim
N_2}\psi\left(\frac{h_2n_2}{h_1}+\frac{\theta(h_1^2-h_2^2)}{2h_1}-\frac{C\Delta
(H_1N_1H_2N_2)^{\frac{1}{4}}}{2h_1}\right),\\
&&\Sigma_3=\sum_{h_1\sim H_1}\sum_{h_2\sim H_2}\sum_{n_2\sim
N_2}\psi\left(\frac{h_2n_2}{h_1}+\frac{\theta(h_1^2-h_2^2)}{2h_1}+\frac{C\Delta
(H_1N_1H_2N_2)^{\frac{1}{4}}}{2h_1}\right).
\end{eqnarray*}
Obviously we have
\begin{equation}
\Sigma_1\ll \Delta (H_1N_1H_2N_2)^{\frac{1}{4}}H_2N_2\ll \Delta
(H_1N_1H_2N_2)^{\frac{3}{4}}
\end{equation}
 if noting $ H_1N_1\asymp H_2N_2.$

Now we estimate $\Sigma_2.$ Taking $H=H_1$ in Lemma 4.1 we get
\begin{eqnarray}
\Sigma_2&&=\sum_{h_1\sim H_1}\sum_{h_2\sim H_2}\sum_{1\leq |h|\leq
H_1}a^{*}(h,h_1,h_2)\sum_{n_2\sim
N_2}e\left(\frac{hh_2n_2}{h_1} \right)\\
&&\ \ \ \ \ \  +O\left( \sum_{h_1\sim H_1}\sum_{h_2\sim
H_2}\sum_{0\leq |h|\leq H_1}b^{*}(h,h_1,h_2)\sum_{n_2\sim
N_2}e\left(\frac{hh_2n_2}{h_1} \right)\right)\nonumber\\
&&\ll H_2N_2\log H_1+\Sigma_4,\nonumber
\end{eqnarray}
where
\begin{eqnarray*}
\Sigma_4&&= \sum_{h_1\sim H_1}\sum_{h_2\sim H_2}\sum_{1\leq h\leq
H_1}\frac{1}{h}\left|\sum_{n_2\sim N_2}e\left(\frac{hh_2n_2}{h_1}
\right)\right|,\\
a^{*}(h,h_1,h_2)&&=a(h)e\left(\frac{\theta(h_1^2-h_2^2)}{2h_1}-\frac{C\Delta
(H_1N_1H_2N_2)^{\frac{1}{4}}}{2h_1}\right)\ll 1/|h|,\\
b^{*}(h,h_1,h_2)&&=b(h)e\left(\frac{\theta(h_1^2-h_2^2)}{2h_1}-\frac{C\Delta
(H_1N_1H_2N_2)^{\frac{1}{4}}}{2h_1}\right)\ll 1/H_1.
\end{eqnarray*}

It is easy to see that
\begin{eqnarray*}
\sum_{n_2\sim N_2}e\left(\frac{hh_2n_2}{h_1}
\right)\ll\left\{\begin{array}{ll}
N_2,&\mbox{$h_1|hh_2,$}\\
\frac{1}{\Vert \frac{hh_2}{h_1}\Vert},& \mbox{$h_1\not|hh_2.$}
\end{array}\right.
\end{eqnarray*}
Thus
\begin{eqnarray}
\Sigma_4&&\ll \Sigma_5+\Sigma_6,\\
\Sigma_5&&=N_2\sum_{h_1\sim H_1}\sum_{1\leq h\leq
H_1}\frac{1}{h}\sum_{\stackrel{ h_2\sim H_2}{hh_2\equiv 0(mod\ h_1)}}1,\nonumber\\
\Sigma_6&&=\sum_{h_1\sim H_1}\sum_{1\leq h\leq
H_1}\frac{1}{h}\sum_{\stackrel{h_2\sim H_2}{hh_2\not\equiv 0(mod\
h_1)}}\frac{1}{\Vert \frac{hh_2}{h_1}\Vert} .\nonumber
\end{eqnarray}

Writing $h=dh^{*}, h_1=dh_1^{*}, d=(h,h_1),$ we get
\begin{eqnarray}
\Sigma_5&&\ll N_2\sum_{d\ll H_1}\frac{1}{d}\sum_{h_1^{*}\sim
\frac{H_1}{d}}\sum_{\stackrel{h^{*}\leq
\frac{H_1}{d}}{(h^{*},h_1^{*})=1}}\frac{1}{h^{*}}\sum_{\stackrel{h_2\sim
H_2}{h_2\equiv 0(mod\ h_1^{*})}}1\\
&&\ll N_2\sum_{d\ll H_1}\frac{1}{d}\sum_{h_1^{*}\sim
\frac{H_1}{d}}\sum_{\stackrel{h^{*}\leq
\frac{H_1}{d}}{(h^{*},h_1^{*})=1}}\frac{1}{h^{*}}\frac{H_2}{h_1^{*}}\nonumber\\
&&\ll N_2H_2\sum_{d\ll H_1}\frac{1}{d} \sum_{ h^{*}\leq
\frac{H_1}{d} }\frac{1}{h^{*}} \nonumber\\
&&\ll N_2H_2\log^2 H_1\ll (H_1N_1H_2N_2)^{1/2}\log^2 H_1\nonumber
\end{eqnarray}
if noting that $H_1N_1\asymp H_2N_2.$

For $\Sigma_6$ we have
\begin{eqnarray}
\Sigma_6&&\ll\sum_{d\ll H_1}\frac{1}{d}\sum_{h_1^{*}\sim
\frac{H_1}{d}}\sum_{\stackrel{h^{*}\leq
\frac{H_1}{d}}{(h^{*},h_1^{*})=1}}\frac{1}{h^{*}}\sum_{\stackrel{h_2\sim
H_2}{h^{*}h_2\not\equiv 0(mod\ h_1^{*})}}\frac{1}{\Vert
\frac{h^{*}h_2}{h_1^{*}}\Vert}.
\end{eqnarray}

We need to bound the sum $$\sum_{\stackrel{h_2\sim
H_2}{h^{*}h_2\not\equiv 0(mod\ h_1^{*})}}\frac{1}{\Vert
\frac{h^{*}h_2}{h_1^{*}}\Vert}.$$

The condition  $h^{*}h_2\not\equiv 0(mod\ h_1^{*})$ implies that
$h_1^{*}\geq 2.$  Let $s=[H_2/h_1^{*}],$ then $sh_1^{*}\leq
H_2<(s+1)h_1^{*}, 2H_2<(2s+2)h_1^{*}.$ Thus
\begin{eqnarray*}
\sum_{\stackrel{h_2\sim H_2}{h^{*}h_2\not\equiv 0(mod\
h_1^{*})}}\frac{1}{\Vert \frac{h^{*}h_2}{h_1^{*}}\Vert} &&\ll
\sum_{j=s}^{2s+1}\sum_{\stackrel{jh_1^{*}<h_2<(j+1)h_1^{*}}{h^{*}h_2\not\equiv
0(mod\ h_1^{*})}}\frac{1}{\Vert \frac{h^{*}h_2}{h_1^{*}}\Vert}\\
&&\ll \sum_{j=s}^{2s+1}\sum_{h_2^{*}=1}^{h_1^{*}-1}\frac{1}{\Vert
\frac{h_2^{*}}{h_1^{*}}\Vert}\ll
\sum_{j=s}^{2s+1}\sum_{h_2^{*}=1}^{[h_1^{*}/2]}\frac{1}{\Vert
\frac{h_2^{*}}{h_1^{*}}\Vert} \\&&\ll
\sum_{j=s}^{2s+1}\sum_{h_2^{*}=1}^{[h_1^{*}/2]}
\frac{h_1^{*}}{h_2^{*}}\ll H_2\log H_1.
\end{eqnarray*}
Inserting this bound into (4.8) we get
\begin{equation}
\Sigma_6\ll H_1H_2\log^2 H_1\ll (H_1N_1H_2N_2)^{1/2}\log^2 H_1
\end{equation}
if noting that $N_1\gg H_1$ and $N_2\gg H_2.$

From (4.4)-(4.9) we get
\begin{equation}
\Sigma_2\ll \Delta (H_1H_2N_1N_2)^{3/4}  +(H_1H_2N_1N_2)^{1/2}\log^2
H_1H_2N_1N_2.
\end{equation}
Similarly we have
\begin{equation}
\Sigma_3\ll \Delta (H_1H_2N_1N_2)^{3/4}  +(H_1H_2N_1N_2)^{1/2}\log^2
H_1H_2N_1N_2.
\end{equation}
Now Lemma 4.4 follows from (4.3), (4.10) and (4.11).
\end{proof}

{\bf Lemma 4.5.} Suppose $\theta>0$ is an irrational number of {\it
approximation type} $\gamma\geq 1.$ Then for any $q\in {\Bbb
N}\bigcup ({\Bbb N}+1/2),$ we have
$$\Vert q\theta\Vert\gg q^{-\gamma-\varepsilon},$$
where the implied constant depending only on $\varepsilon.$

\begin{proof} It follows easily from the definition of
{\it approximation type}.

\end{proof}

{\bf Lemma 4.6.} Suppose $\theta>0$ is an irrational number of {\it
approximation type} $\gamma\geq 1$ and that
\begin{eqnarray}
&&0<|\alpha(\theta;h_1,h_2,n_1,n_2)|<\frac{1}{10}\beta^{1/4}(\theta;h_1,h_2,n_1,n_2),
\ n_1>\theta h_1, n_2>\theta h_2,
\end{eqnarray}
where
$$\beta(\theta;h_1,h_2,n_1,n_2):=
 h_1(2n_1-\theta h_1)h_2(2n_2-\theta h_2).$$

If $h_1=h_2=h, $ then
$$|\alpha(\theta;h ,h ,n_1,n_2)|\gg \frac{h^{1/2} }{(n_1n_2)^{1/4}}.$$

If $h_1\not= h_2,$ then
$$|\alpha(\theta;h_1,h_2,n_1,n_2)|\gg \frac{|h_1h_2|^{\frac{-\gamma-\varepsilon}{2}}}{(h_1h_2n_1n_2)^{1/4}}.$$

\begin{proof}
The condition (4.12)  implies that
$$h_1(2n_1-\theta
h_1)\asymp h_2(2n_2-\theta h_2)\asymp h_1n_1\asymp h_2n_2.$$ So we
have
\begin{eqnarray}
|\alpha(\theta;h_1,h_2,n_1,n_2|&&=\frac{|2(h_1n_1-h_2n_2)-\theta
(h_1^2-h_2^2)|}{\sqrt{h_1(2n_1-\theta h_1)}+\sqrt{h_2(2n_2-\theta
h_2)}}\\
&&\gg \frac{|2(h_1n_1-h_2n_2)-\theta
(h_1^2-h_2^2)|}{(h_1h_2n_1n_2)^{1/4}}\nonumber
\end{eqnarray}

If $h_1=h_2=h,$ then (4.13) implies
\begin{eqnarray*}
|\alpha(\theta;h_1,h_2,n_1,n_2| \gg
\frac{h|n_1-n_2|}{h^{1/2}(n_1n_2)^{1/4}}\gg \frac{h^{1/2} }{
(n_1n_2)^{1/4}}.
\end{eqnarray*}

If $h_1\not= h_2,$ then from (4.13) and Lemma 4.5  we have
\begin{eqnarray}
|\alpha(\theta;h_1,h_2,n_1,n_2|&& \gg \frac{\Vert\theta
(h_1^2-h_2^2)\Vert}{(h_1h_2n_1n_2)^{1/4}}\nonumber\\
&&\gg
\frac{|h_1^2-h_2^2|^{-\gamma-\varepsilon}}{(h_1h_2r_1r_2)^{1/4}}\nonumber\\
&&\gg
\frac{|h_1+h_2|^{-\gamma-\varepsilon}|h_1-h_2|^{-\gamma-\varepsilon}}{(h_1h_2r_1r_2)^{1/4}}\nonumber\\
&&\gg \frac{|h_1+h_2|^{-\gamma-\varepsilon}
}{(h_1h_2r_1r_2)^{1/4}}\nonumber\\
&&\gg \frac{|h_1h_2|^{\frac{-\gamma-\varepsilon}{2}}
}{(h_1h_2r_1r_2)^{1/4}},\nonumber
\end{eqnarray}
where in the last step we used the inequality $a^2+b^2\geq 2ab.$
\end{proof}

\section{\bf Proof of Theorem 1 }

In this section we shall prove Theorem 1.   It suffices for us to
evaluate the integral $\int_T^{2T}R^2(2\pi x)dx$ , where $T$ is a
large real number.

Suppose $H:=T^2.$  By Lemma 3.1 and Lemma  4.2 we have
\begin{eqnarray}
&&\ \ \ \ \ \ \ R(2\pi x)=R_1(x,H)+R_2(x,H),\\
&&R_1(x;H)=\frac{2^{1-l}}{(l-1)!\pi i}\sum_{1\leq |h|\leq
H}\frac{e(-lh/2)}{h}\sum_{m\leq
\sqrt \frac{x}{\theta}}m(x-\theta m^2)^{l-1}e\left(h\left(\frac{x}{2m}-\frac{\theta m}{2}\right)\right),\nonumber\\
&&R_2(x,H)=O(T^{l-1/2}G(x,H)+T^{l-1/2}),\nonumber\\ &&G(x,H)=
\sum_{m\leq \sqrt{\frac{x}{\theta}}}\min\left(1,\frac{1}{H\Vert
\frac{x}{2m}-\frac{\theta m}{2}+\frac{l}{2}\Vert}\right).\nonumber
\end{eqnarray}

\subsection{\bf The mean-square  of  $R_1(x,H)$}\

In this subsection we shall study The mean-square  of $R_1(x,H)$.
First we shall derive an Vorono\"{i} type formula for $R_1(x,H).$

\subsubsection{\bf Vorono\"{i} analogue formula for $R_1(x,H)$}

We have
 \begin{eqnarray}
\ \ \ \ R_1(x,H)&&=\frac{2^{1-l}}{(l-1)!\pi i}\sum_{1\leq |h|\leq
H}\frac {e(-lh/2)}{h}\sum_{m\leq \sqrt\frac{x}{\theta}}m(x-\theta
m^2)^{l-1}\\&&\ \ \ \ \ \ \ \ \ \ \ \ \ \ \  \ \  \times
e\left(h\left(\frac{x}{2m}-\frac
{\theta m}{2}\right)\right)\nonumber\\
&&=\frac{2^{1-l}}{(l-1)!\pi
i}\sum_{j_1=0}^{l-1}(-\theta)^{j_1}{l-1\choose j_1}\sum_{1\leq
|h|\leq
H}\frac {e(-lh/2)}{h}\sum_{m\leq \sqrt{\frac{x}{\theta}}}x^{l-1-j_1}\nonumber\\
&&\ \ \ \ \ \ \ \ \ \ \ \ \ \ \ \ \ \ \ \ \ \ \ \ \ \ \  \times
m^{2j_1+1}e\left(h\left(\frac{x}{2m}-\frac{\theta m}{2} \right)\right)\nonumber\\
&&=\frac{2^{1-l}}{(l-1)!}\sum_{j_1=0}^{l-1}(-\theta)^{j_1}{l-1\choose
j_1}F(x;j_1),\nonumber
 \end{eqnarray}
say, where
$$F(x;j_1):=\frac{1}{\pi
i} \sum_{1\leq |h|\leq H}\frac {e(-lh/2)}h\sum_{m\leq \sqrt
\frac{x}{\theta}}x^{l-1-j_1}m^{2j_1+1}
e\left(h\left(\frac{x}{2m}-\frac {\theta m}{2} \right)\right).$$

Let $J=[({\cal L}-\log {\cal L})/2\log 2]$   we get
\begin{eqnarray}
&&\ \ \ \ \ \ \ \ \ \ \ F(x;j_1)\\&&=\frac{1}{\pi i}\sum_{-H\leq
h\leq -1}\frac{ e(-lh/2)}{h}\sum_{j=0}^J\sum_{m\sim \sqrt
\frac{x}{\theta}2^{-j-1}}x^{l-1-j_1}m^{2j_1+1}e\left(h\left(\frac{x}{2m}-\frac{\theta m}{2}\right)\right)\nonumber \\
&&\hspace{4mm}+\frac{1}{\pi i}\sum_{1\leq h\leq
H}\frac{e(-lh/2)}{h}\sum_{j=0}^J\sum_{m\sim \sqrt
\frac{x}{\theta}2^{-j-1}}x^{l-1-j_1}m^{2j_1+1}e\left(h\left(\frac{x}{2m}-\frac{\theta
m}{2}\right)\right)+O(x^{l-1}{\cal
L}^2)\nonumber\\
&&= -\frac{1}{\pi i}\sum_{1\leq h\leq H}\frac{
e(lh/2)}{h}\sum_{j=0}^J\sum_{m\sim \sqrt
\frac{x}{\theta}2^{-j-1}}x^{l-1-j_1}m^{2j_1+1}e\left(-h\left(\frac{x}{2m}-\frac{\theta m}{2}\right)\right)\nonumber\\
&&\hspace{4mm} +\frac{1}{\pi i}\sum_{1\leq h\leq
H}\frac{e(-lh/2)}{h}\sum_{j=0}^J\sum_{m\sim \sqrt
\frac{x}{\theta}2^{-j-1}}x^{l-1-j_1}m^{2j_1+1}e\left(h\left(\frac{x}{2m}-\frac{\theta
m}{2}\right)\right)+O(x^{l-1}{\cal
L}^2) \nonumber\\
&&=-\frac{\Sigma_7}{\pi i}+\frac{\overline{\Sigma_7}}{\pi i}
+O(x^{l-1}{\cal L}^2),\nonumber
\end{eqnarray}
where
\begin{eqnarray*}
\Sigma_7=\sum_{1\leq h\leq H}\frac{
e(lh/2)}{h}\sum_{j=0}^J\sum_{m\sim \sqrt
\frac{x}{\theta}2^{-j-1}}x^{l-1-j_1}m^{2j_1+1}e\left(-h\left(\frac{x}{2m}-\frac{\theta
m}{2}\right)\right).
\end{eqnarray*}

Let $$  S(x;h,j_1,j)=\sum_{m\sim \sqrt
\frac{x}{\theta}2^{-j-1}}x^{l-1-j_1}m^{2j_1+1}e\left(-h\left(\frac{x}{2m}-\frac{\theta
m}{2}\right)\right).$$

 By Lemma 4.3 we get
\begin{eqnarray}
 \ \ \ S(x;h,j_1,j)&&=e^{-\frac{\pi i}{4}}\sum_{\beta(h,j)<r\leq
\beta(h,j+1)}\frac{x^{l-1/4}h^{j_1+3/4}}{(2r-\theta
h)^{j_1+5/4}}e\left(-\sqrt{xh(2r-\theta h)}\right)\\
&&+O(x^{l-1/2}{\cal L})+O\left(x^{l-1/2}\min\left(\frac{x^{\frac
14}}{h^{\frac 12}2^{\frac{3j}{2}}},\frac{1}{\Vert \beta(h,j)
\Vert}\right)\right)\nonumber\\
&&+O\left(x^{l-1/2}\min\left(\frac{x^{\frac 14}}{h^{\frac
12}2^{\frac{3j}{2}}},\frac{1}{\Vert \beta(h,j+1)
\Vert}\right)\right),\nonumber
\end{eqnarray}
where
\begin{equation}
\beta(h,j): =\theta h(2^{2j-1}+1/2).
\end{equation}

Inserting (5.4) into $\Sigma_7$ we have
\begin{eqnarray}
 \ \ \ \ \ \   \ \  \Sigma_7&&=e^{-\pi i/4}\sum_{1\leq h\leq
H}\frac{e(lh/2)}{h}\sum_{j=0}^J\sum_{\beta(h,j)<r\leq
 \beta(h,j+1)}\\
&&\hspace{15mm}\frac{x^{l-1/4}h^{j_1+3/4}}{(2r-\theta
h)^{j_1+5/4}}e(-\sqrt{xh(2r-\theta h)})+O( g(T,H)) \nonumber\\&&=
\sum_{1\leq h\leq H}\frac{e(lh/2)}{h} \sum_{\theta h<r\leq
\theta h(2^{2J+1}+1/2)}\frac{x^{l-1/4}h^{j_1+3/4}}{(2r-\theta h)^{j_1+5/4}}e(-\sqrt{xh(2r-\theta h)}-\frac{1}{8})\nonumber\\
&&\ \  \ \ \ \ +O(g(T,H)),\nonumber
\end{eqnarray}
where
\begin{eqnarray*}
g(T,H): = T^{l-1/2}{\cal L}^3+T^{l-1/2}\sum_{1\leq h\leq
H}\frac{1}{h}\sum_{j=0}^J \min\left(\frac{T^{\frac 14}}{h^{\frac
12}2^{\frac{3j}{2}}},\frac{1}{\Vert \beta(h,j) \Vert}\right).
\end{eqnarray*}
By Lemma 4.5 we have $\Vert \beta(h,j) \Vert\gg
(h2^{2j})^{-\gamma-\varepsilon},$ which implies that

\begin{eqnarray}
\ \ g(T,H)&&\ll T^{l-1/2}{\cal L}^3+T^{l-1/2}\sum_{1\leq h\leq
H}\frac{1}{h}\sum_{j=0}^J \min\left(\frac{T^{\frac 14}}{h^{\frac
12}2^{\frac{3j}{2}}}, (h2^{2j})^{\gamma+\varepsilon} \right)\\
&&\ll T^{l-1/2}{\cal L}^3+T^{l-1/2}\sum_{1\leq h\leq
H}\frac{1}{h}\sum_{j=0}^J\left(\frac{T^{\frac 14}}{h^{\frac
12}2^{\frac{3j}{2}}}\right)^{\frac{2\gamma}{2\gamma+1}}
\left((h2^{2j})^{\gamma+\varepsilon}\right)^{\frac{1}{2\gamma+1}}\nonumber\\
&&\ll T^{l-1/2+\frac{\gamma}{4\gamma+2}+\varepsilon}.\nonumber
\end{eqnarray}

Inserting (5.4)-(5.7) into (5.3) we get
\begin{eqnarray}
&& \ \ \ \ \ \ \ \ \ \ F(x;j_1)=\frac{1}{\pi i}\sum_{1\leq h\leq
H}\frac{e(lh/2)}{h}   \sum_{\theta h<r\leq
\theta h(2^{2J+1}+1/2)}\frac{x^{l-1/4}h^{j_1+3/4}}{(2r-\theta h)^{j_1+5/4}}\\
&&\ \ \ \ \ \ \ \ \ \ \ \ \  \times \left(e(\sqrt{xh(2r-\theta
h)}+\frac
18)-e(\sqrt{xh(2r-\theta h)}-\frac 18)\right)+O(T^{l-1/2+\frac{\gamma}{4\gamma+2}+\varepsilon})\nonumber\\
&&=\frac{2x^{l-1/4}}{\pi }\sum_{1\leq h\leq H}\frac{e(lh/2)}{h}
\sum_{\theta h<r\leq \theta h(2^{2J+1}+1/2)}\frac{
h^{j_1+3/4}}{(2r-\theta h)^{j_1+5/4}}\sin(2\pi\sqrt{xh(2r-\theta h)}+\frac{\pi}{4})\nonumber\\
&&\ \ \ \ \ \ \ \ \ \ \ \ \ \ \ +O(T^{l-1/2+\frac{\gamma}{4\gamma+2}+\varepsilon})\nonumber\\
&&=\frac{2x^{l-1/4}}{\pi }\sum_{1\leq h\leq H}\frac{e(lh/2)}{h}
\sum_{\theta h<r\leq \theta h(2^{2J+1}+1/2)}\frac{
h^{j_1+3/4}}{(2r-\theta h)^{j_1+5/4}}\cos(2\pi\sqrt{xh(2r-\theta h)}-\frac{\pi}{4})\nonumber\\
&&\ \ \ \ \ \ \ \ \ \ \ \ \ \ \
+O(T^{l-1/2+\frac{\gamma}{4\gamma+2}+\varepsilon}).\nonumber
\end{eqnarray}

From (5.2) and (5.8) we get
\begin{eqnarray}
&&
 R_1(x,H)=\frac{2^{2-l}x^{l-1/4}}{(l-1)!\pi}\sum_{1\leq h\leq H}
   \sum_{\theta h<r\leq
\theta h(2^{2J+1}+1/2)}\cos(2\pi\sqrt{xh(2r-\theta h)}-\frac{\pi}{4}) \\
&& \ \ \ \ \ \ \  \ \ \ \ \ \ \ \ \ \ \times
\sum_{j_1=0}^{l-1}(-1)^{j_1}{l-1\choose j_1}\frac{
e(lh/2)h^{j_1-1/4}}{(2r-\theta h)^{j_1+5/4}}+O( T^{l-1/2+\frac{\gamma}{4\gamma+2}+\varepsilon} )\nonumber\\
&&\ \ \ \ \ \ \ \ \ \ \ \ =
R_{11}(x,H)+O(T^{l-1/2+\frac{\gamma}{4\gamma+2}+\varepsilon})
,\nonumber
\end{eqnarray}
where
\begin{eqnarray}
 R_{11}(x,H):&&=\frac{2^{2-l}x^{l-1/4}}{(l-1)!\pi}\sum_{1\leq h\leq H}
  \sum_{\theta h<r\leq
\theta h(2^{2J+1}+1/2)}u(h,r)\\
&&\ \ \ \ \ \ \ \ \ \ \ \ \times \cos(2\pi\sqrt{xh(2r-\theta
h)}-\frac{\pi}{4}),\nonumber\\ u(h,r):&& =
\frac{e(lh/2)}{h^{1/4}(2r-\theta h)^{5/4}}\left(1-\frac{\theta
h}{2r-\theta h}\right)^{l-1}.\nonumber
\end{eqnarray}

\subsubsection{\bf Mean square of $R_{11}(x,H)$}

Now we study the mean square of $R_{11}(x,H)$. By the elementary
formula
$$\cos u\cos v=\frac{1}{2}(\cos(u-v)+\cos(u+v))$$
we have
\begin{eqnarray}
R^2_{11}(x,H)&&=\frac{2^{4-2l}x^{2l-\frac{1}{2}}}{\pi^2(l-1)!^2}\sum_{8}u(h_1,r_1)u(h_2,r_2)\\
&&\ \ \ \ \ \ \ \
\times\prod_{j=1}^2\cos\left(2\pi\sqrt{xh_j(2r_j-\theta
h_j)}-\frac{\pi}{4}\right)\nonumber\\
&&=S_1(x)+S_2(x)+S_3(x),\nonumber
\end{eqnarray}
where
\begin{eqnarray*}
&&S_1(x)=\frac{2^{3-2l}x^{2l-\frac{1}{2}}}{\pi^2(l-1)!^2}\sum_{9}u(h_1,r_1)u(h_2,r_2),\\
&&S_2(x)=\frac{2^{3-2l}x^{2l-\frac{1}{2}}}{\pi^2(l-1)!^2}\sum_{10}u(h_1,r_1)u(h_2,r_2)
\cos(\alpha(\theta;h_1,h_2,r_1,r_2)),\\
&&S_3(x)=\frac{2^{3-2l}x^{2l-\frac{1}{2}}}{\pi^2(l-1)!^2}\sum_{8}u(h_1,r_1)u(h_2,r_2)\\
&&\ \ \ \ \ \ \ \times \sin\left(2\pi\sqrt{xh_1(2r_1-\theta
h_1)}+2\pi\sqrt{xh_2(2r_2-\theta h_2)}\right),\\
&&SC(\Sigma_8): 1\leq h_j\leq H, \theta h_j<r_j\leq \theta
h_j(2^{2J+1}+1/2)\ \ (j=1,2), \\
&&SC(\Sigma_9): 1\leq h_j\leq H, \theta h_j<r_j\leq \theta
h_j(2^{2J+1}+1/2)\ \ (j=1,2), \\
&&\hspace{20mm} \alpha(\theta;h_1,h_2,r_1,r_2)=0,\\
&&SC(\Sigma_{10}): 1\leq h_j\leq H, \theta h_j<r_j\leq \theta
h_j(2^{2J+1}+1/2)\ \ (j=1,2),\\
&&\hspace{20mm} \alpha(\theta;h_1,h_2,r_1,r_2)\not= 0.
\end{eqnarray*}

We first consider the contribution of $S_1(x).$  Since $\theta $ is
irrational ,  we see easily that $\alpha(\theta;h_1,h_2,r_1,r_2)=0$
holds if and only if $h_1=h_2, r_1=r_2.$ Thus we have
\begin{eqnarray}
\ \ \ \ \ \int_T^{2T}S_1(x)dx
=\frac{2^{3-2l}}{\pi^2(l-1)!^2}\sum_{1\leq h\leq H}\sum_{\theta
h<r\leq \theta
h(2^{2J+1}+1/2)}u^2(h,r)\int_T^{2T}x^{2l-\frac{1}{2}}dx.
\end{eqnarray}

Recalling the definition of $u(h,r)$ , we have
\begin{eqnarray}
&&\ \ \ \ \ \ \sum_{1\leq h\leq H}\sum_{\theta h<r\leq \theta
h(2^{2J+1}+1/2)}u^2(h,r)\\
&&=\sum_{1\leq h\leq H}\sum_{\theta h<r }\frac{(1-\frac{\theta
h}{2r-\theta h})^{2l-2}}{h^{1/2}(2r-\theta
h)^{5/2}}+O\left(\sum_{1\leq h\leq H}\sum_{r>\theta
h2^{2J+1}}h^{-1/2}r^{-5/2}\right)\nonumber\\
&&=\sum_{h=1}^\infty\sum_{r>\theta h }\frac{(1-\frac{\theta
h}{2r-\theta h})^{2l-2}}{h^{1/2}(2r-\theta
h)^{5/2}}+O(T^{-3/2}+H^{-2})\nonumber\\
&&={\cal C}_{l,\theta}+O(T^{-3/2} ),\nonumber
\end{eqnarray}
which combining (5.12) gives
\begin{eqnarray}
\int_T^{2T}S_1(x)dx=\frac{2^{3-2l}{\cal
C}_{l,\theta}}{\pi^2(l-1)!^2}\int_T^{2T}x^{2l-\frac{1}{2}}dx
+O(T^{2l-1} ).
\end{eqnarray}

For the contribution of $S_3(x),$ by the first derivative test we
get
\begin{eqnarray}
&& \ \ \ \ \ \int_T^{2T}S_3(x)dx\ll
T^{2l}\sum_{8}\frac{(\sqrt{h_1(2r_1-\theta
h_1)}+\sqrt{h_2(2r_2-\theta
h_2)})^{-1}}{(h_1h_2)^{\frac{1}{4}}(2r_1-\theta h_1)^{\frac
54}(2r_2-\theta h_2)^{\frac{5}{4}}}\\
&&\ll T^{2l}\sum_{8}\frac{(\sqrt{h_1(2r_1-\theta
h_1)}\sqrt{h_2(2r_2-\theta
h_2)})^{-1/2}}{(h_1h_2)^{\frac{1}{4}}(2r_1-\theta h_1)^{\frac
54}(2r_2-\theta h_2)^{\frac{5}{4}}}\nonumber\\
&&\ll T^{2l}\sum_{8}\frac{1}{(h_1h_2)^{\frac{1}{2}}(2r_1-\theta
h_1)^{\frac 32}(2r_2-\theta h_2)^{\frac{3}{2}}}\nonumber\\
&&\ll T^{2l}\left(\sum_{1\leq h\leq H}h^{-1/2}\sum_{\theta h<r\leq
\theta h(2^{2J+1}+1/2)}r^{-3/2}\right)^2\nonumber\\
&&\ll T^{2l}{\cal L}^2,\nonumber
\end{eqnarray}
where in the second step we used the inequality $a^2+b^2\geq 2ab.$

Finally we consider the contribution of $S_2(x).$ By the first
derivative test we get
\begin{eqnarray}
 &&\ \ \ \ \ \ \ \ \ \ \ \ \int_T^{2T}S_2(x)dx\\&&\ll
T^{2l-1/2}\sum_{10}\frac{1}{h_1^{1/4}h_2^{1/4}r_1^{5/4}r_2^{5/4}}\min\left(T,\frac{T^{1/2}}{|\alpha(\theta;h_1,h_2,r_1,r_2)|}\right)\nonumber\\
&&\ll
T^{2l}\sum_{11}\frac{1}{h_1^{1/4}h_2^{1/4}r_1^{5/4}r_2^{5/4}}\min\left(T^{1/2},\frac{ 1}{|\alpha(\theta;h_1,h_2,r_1,r_2)|}\right)\nonumber\\
&&\ \ \ +T^{2l}\sum_{12}\frac{1}{h^{1/2}r_1^{5/4}r_2^{5/4}}\min\left(T^{1/2},\frac{ 1}{|\alpha(\theta;h,h,r_1,r_2)|}\right)\nonumber\\
&&\ \ \
+T^{2l}\sum_{13}\frac{1}{h_1^{1/4}h_2^{1/4}r_1^{5/4}r_2^{5/4}}\min\left(T^{1/2},\frac{
1}{|\alpha(\theta;h_1,h_2,r_1,r_2)|}\right),\nonumber
\end{eqnarray}
where
\begin{eqnarray*}
 &&SC(\Sigma_{11}): 1\leq h_j\leq H, \theta h_j<r_j\leq \theta
h_j(2^{2J+1}+1/2)\ \ (j=1,2),\\
&&\hspace{17mm} |\alpha(\theta;h_1,h_2,r_1,r_2)|\geq
\frac{1}{10} \beta^{1/4}(\theta;h_1,h_2,r_1,r_2),\\
&&SC(\Sigma_{12}): 1\leq h\leq H, \theta h<r_j\leq \theta
h(2^{2J+1}+1/2)\ (j=1,2), r_1\not= r_2,\\
&&\hspace{17mm} |\alpha(\theta;h,h,r_1,r_2)|<
\frac{1}{10} \beta^{1/4}(\theta;h,h,r_1,r_2),\\
&&SC(\Sigma_{13}): 1\leq h_j\leq H, \theta h_j<r_j\leq \theta
h_j(2^{2J+1}+1/2)\ \ (j=1,2),h_1\not= h_2,\\
&&\hspace{17mm} |\alpha(\theta;h_1,h_2,r_1,r_2)|<
\frac{1}{10}\beta^{1/4}(\theta;h_1,h_2,r_1,r_2)
\end{eqnarray*}
and where $\beta(\theta;h_1,h_2,r_1,r_2)$ was defined in Lemma 4.6.

Similar to the case $S_3(x)$, we have
\begin{eqnarray}
&&\ \ \ \ \
 T^{2l}\sum_{11}\frac{1}{h_1^{1/4}h_2^{1/4}r_1^{5/4}r_2^{5/4}}\min\left(T^{1/2},\frac{1}{|\alpha(\theta;h_1,h_2,r_1,r_2)|}\right)
 \\
&&\ll
T^{2l}\sum_{11}\frac{1}{h_1^{1/4}h_2^{1/4}r_1^{5/4}r_2^{5/4}}\times\frac{1}{|\alpha(\theta;h_1,h_2,r_1,r_2)|}\nonumber\\
&&\ll T^{2l}\sum_{11}\frac{1}{h_1^{1/2}h_2^{1/2}r_1^{3/2}r_2^{3/2}}
\ll T^{2l}{\cal L}^{2}\nonumber.
\end{eqnarray}

By Lemma 4.6 we have $|\alpha(\theta;h,h,r_1,r_2)|\gg
h^{1/2}(r_1r_2)^{-1/4}$ under the condition $SC(\Sigma_{12}).$ Thus
\begin{eqnarray}
&&\ \ \ \ \ \ \
 T^{2l}\sum_{12}\frac{1}{h^{1/2}r_1^{5/4}r_2^{5/4}}\min\left(T^{1/2},\frac{1}{|\alpha(\theta;h,h,r_1,r_2)|}\right)
 \\
&&\ll
T^{2l}\sum_{12}\frac{1}{h^{1/2}r_1^{5/4}r_2^{5/4}}\times\frac{1}{|\alpha(\theta;h,h,r_1,r_2)|}\nonumber\\
&&\ll T^{2l}\sum_{h,r_1,r_2}\frac{1}{hr_1r_2} \ll T^{2l}{\cal
L}^{3}\nonumber.
\end{eqnarray}

Now we estimate the sum
 $\Sigma_{13}.$ By a splitting argument we have
\begin{eqnarray}
&&\ \ \ \ \ \ \ \ \sum_{13}\frac{1}{h_1^{1/4}h_2^{1/4}r_1^{5/4}r_2^{5/4}}\min\left(T^{1/2},\frac{1}{|\alpha(\theta;h_1,h_2,r_1,r_2)|}\right)\\
&&\ll {\cal
L}^4\sum_{14}\frac{1}{h_1^{1/4}h_2^{1/4}r_1^{5/4}r_2^{5/4}}\min\left(T^{1/2},\frac{1}{|\alpha(\theta;h_1,h_2,r_1,r_2)|}\right)\nonumber\\
&&\ll (H_1H_2)^{-1/4}(N_1N_2)^{-5/4}{\cal L}^4\sum_{14}
\min\left(T^{1/2},\frac{1}{|\alpha(\theta;h_1,h_2,r_1,r_2)|}\right)\nonumber\\
&&\ll U_1+U_2\nonumber
\end{eqnarray}
for some $(H_1,H_2,N_1,N_2)$  for which $$1\ll H_j\ll H, H_j\ll
N_j\ll H2^{2J}(j=1,2), $$ where
\begin{eqnarray*}
&&U_1=T^{1/2}(H_1H_2)^{-1/4}(N_1N_2)^{-5/4}{\cal L}^4\times {\cal A}_\theta(H_1,H_2,N_1,N_2;T^{-1/2}),\\
&&U_2=(H_1H_2)^{-1/4}(N_1N_2)^{-5/4}{\cal
L}^4\sum_{15}\frac{1}{|\alpha(\theta;h_1,h_2,r_1,r_2)|},
\end{eqnarray*}
\begin{eqnarray*}
&& SC(\Sigma_{14}):  h_j\sim H_j, r_j\sim N_j(j=1,2),  \ h_1\not= h_2,\\
&&\ \ \ \ \ \ \ \ \ \ \ 0<|\alpha(\theta;h_1,h_2,r_1,r_2)|<
\frac{1}{10} \beta^{1/4}(\theta;h_1,h_2,r_1,r_2),\\
&& SC(\Sigma_{15}) :   h_j\sim H_j, r_j\sim N_j(j=1,2),  \ h_1\not= h_2,\\
&&\ \ \ \ \ \ \ \ \ \ \
T^{-1/2}<|\alpha(\theta;h_1,h_2,r_1,r_2)|<\frac{1}{10}\beta^{1/4}(\theta;h_1,h_2,r_1,r_2).
\end{eqnarray*}
We first estimate $U_1.$ By Lemma 4.4 we have
\begin{eqnarray*}
U_1&&\ll T^{1/2}(H_1H_2)^{-1/4}(N_1N_2)^{-5/4}{\cal
L}^4\left(T^{-1/2}(H_1H_2N_1N_2)^{3/4}+(H_1H_2N_1N_2)^{1/2}{\cal
L}^2\right)\\
&&\ll (H_1H_2)^{1/2}(N_1N_2)^{-1/2}{\cal
L}^4+T^{1/2}(H_1H_2)^{1/4}(N_1N_2)^{-3/4}{\cal L}^6\\
&&\ll {\cal L}^4+T^{1/2} (N_1N_2)^{-1/2}{\cal L}^6
\end{eqnarray*}
if noting that $N_j\gg H_j(j=1,2).$ Suppose $(h_1,h_2,r_1,r_2)$
satisfies the conditions of $\Sigma_{14}.$ By Lemma 4.6 we have
\begin{eqnarray*}
|\alpha(\theta;h_1,h_2,r_1,r_2)| \gg (H_1H_2)^{-\frac
{1}{4}-\frac{\gamma}{2}-\frac{\varepsilon}{2}}(N_1N_2)^{-\frac{1}{4}},
\end{eqnarray*}
which combining $|\alpha(\theta;h_1,h_2,r_1,r_2)|\leq T^{-1/2}$
gives
\begin{eqnarray*}
(H_1H_2)^{ \frac{\gamma}{2}+\frac
14+\frac{\varepsilon}{2}}(N_1N_2)^{\frac{1}{4}}\gg T^{1/2}.
\end{eqnarray*}
Hence
\begin{eqnarray*}
(N_1N_2)^{\frac 12+ \frac{\gamma}{2} +\frac{\varepsilon}{2}}
\gg(H_1H_2)^{ \frac{\gamma}{2}+\frac
14+\frac{\varepsilon}{2}}(N_1N_2)^{\frac{1}{4}}\gg T^{1/2},
\end{eqnarray*}
namely
$$N_1N_2\gg T^{\frac{1}{1+\gamma+\varepsilon}}.$$
From the above estimates we get
\begin{equation}
U_1\ll T^{1/2-1/2(1+\gamma)+\varepsilon}.
\end{equation}

Now we estimate $U_2.$ By a splitting argument we have
\begin{eqnarray*}
 U_2\ll \eta^{-1}(H_1H_2)^{-1/4}(N_1N_2)^{-5/4}{\cal
L}^5 \times {\cal A}_\theta(H_1,H_2,N_1,N_2;\eta)
\end{eqnarray*}
for some $T^{-1/2}\ll \eta\ll (H_1H_2N_1N_2)^{1/4}.$  By Lemma 4.6
we get
\begin{eqnarray*}
 U_2&&\ll \eta^{-1}(H_1H_2)^{-1/4}(N_1N_2)^{-5/4}{\cal
L}^5 (\eta (H_1H_2N_1N_2)^{3/4}+(H_1H_2N_1N_2)^{1/2}{\cal L}^2)\\
&&\ll (H_1H_2)^{1/2}(N_1N_2)^{-1/2}{\cal L}^5+\eta^{-1}
(H_1H_2)^{1/4}(N_1N_2)^{-3/4}{\cal L}^7\\
&&\ll {\cal L}^5+\eta^{-1} (H_1H_2)^{1/4}(N_1N_2)^{-3/4}{\cal L}^7.
\end{eqnarray*}

From Lemma 4.6 we get
\begin{eqnarray*}
\eta\gg (H_1H_2)^{-\frac
14-\frac{\gamma}{2}-\frac{\varepsilon}{2}}(N_1N_2)^{-\frac{1}{4}},
\end{eqnarray*}
which combining $\eta\gg T^{-1/2}$ gives
$$\eta^{-1}\ll \min\left(T^{1/2},(H_1H_2)^{\frac{\gamma}{2}+\frac
14+\frac{\varepsilon}{2}}(N_1N_2)^{\frac{1}{4}}\right).$$ Thus we
have
\begin{eqnarray}
 U_2&&\ll  {\cal L}^5+  \min\left(T^{1/2},(H_1H_2)^{\frac{\gamma}{2}+\frac
14+\frac{\varepsilon}{2}}(N_1N_2)^{\frac{1}{4}}\right)\\&&\ \ \
 \ \ \ \ \ \ \ \ \  \ \ \ \ \ \ \times
(H_1H_2)^{1/4}(N_1N_2)^{-3/4}{\cal L}^7\nonumber\\
&&\ll {\cal
L}^5+\min\left(T^{1/2}(H_1H_2)^{1/4}(N_1N_2)^{-3/4},(H_1H_2)^{\frac{1+\gamma}{2}
+\frac{\varepsilon}{2}}(N_1N_2)^{-\frac{1}{2}}\right)\nonumber\\
&&\ll {\cal L}^5+\min\left(T^{1/2}(H_1H_2)^{-1/2}
,(H_1H_2)^{\frac{\gamma}{2}  +\frac{\varepsilon}{2}}
\right)\nonumber\\
&&\ll {\cal
L}^5+(T^{1/2}(H_1H_2)^{-1/2})^{\gamma/(1+\gamma)}((H_1H_2)^{\frac{\gamma}{2} +\frac{\varepsilon}{2}})^{1/(1+\gamma)}\nonumber\\
&&\ll T^{1/2-1/2(1+\gamma)+\varepsilon}.\nonumber
\end{eqnarray}

From (5.16)-(5.21) we get
\begin{equation}
\int_T^{2T}S_2(x)dx\ll T^{2l+1/2-1/2(1+\gamma)+\varepsilon},
\end{equation}
which combining (5.11), (5.14) and (5.15) implies that
\begin{eqnarray}
\int_T^{2T}R_{11}^2(x,H)dx&&=\frac{2^{3-2l}C_{l,\theta}}{\pi^2(l-1)!^2}\int_T^{2T}x^{2l-\frac{1}{2}}dx
 +O(T^{2l+1/2-1/2(1+\gamma)+\varepsilon} ).
\end{eqnarray}

\subsubsection{\bf Mean square of $R_1(x,H)$}

We have
\begin{eqnarray}
\ \ \ \ \ \
R_1^2(x,H)=R_{11}^2(x,H)+O(|R_{11}(x,H)|T^{l-1/2+\frac{\gamma}{4\gamma+2}+\varepsilon}+T^{2l-1+\frac{2\gamma}{4\gamma+2}+2\varepsilon}).
\end{eqnarray}
By (5.7) , (5.23) and Cauchy's inequality we get
\begin{equation}
T^{l-1/2+\frac{\gamma}{4\gamma+2}}\int_T^{2T}|R_{11}(x,H)|dx\ll
T^{2l +\frac{4\gamma+1}{8\gamma+4}+\varepsilon}.
\end{equation}

From (5.23), (5.24) and (5.25) we get
\begin{eqnarray}
\int_T^{2T}R_{1}^2(x,H)dx&&=\frac{2^{3-2l}C_{l,\theta}}{\pi^2(l-1)!^2}\int_T^{2T}x^{2l-\frac{1}{2}}dx
 +O(T^{2l +\frac{4\gamma+1}{8\gamma+4} +\varepsilon} ).
\end{eqnarray}

\subsection{\bf Mean square of $R_2(x,H)$}

 We first study the integral $\int_T^{2T}G(x,H)dx.$
 We have
 \begin{eqnarray}
\int_T^{2T}G(x,H)dx&&\ll \int_T^{2T}\sum_{m\leq
\sqrt{\frac{2T}{\theta}}}\min\left(1,\frac{1}{H\Vert
\frac{x}{2m}-\frac{\theta m}{2}+\frac l2 \Vert}\right)dx\\
&&\ll \sum_{m\leq
\sqrt{\frac{2T}{\theta}}}\int_T^{2T}\min\left(1,\frac{1}{H\Vert
\frac{x}{2m}-\frac{\theta m}{2}+\frac l2 \Vert}\right)dx\nonumber\\
 &&\ll \sum_{m\leq
\sqrt{\frac{2T}{\theta}}}m\int_{\frac{T}{2m}-\frac{\theta
m}{2}+\frac l2 }^{\frac{T}{m}-\frac{\theta m}{2}+\frac l2
}\min\left(1,\frac{1}{H\Vert u  \Vert}\right)du\nonumber\\
&&\ll \sum_{m\leq \sqrt{\frac{2T}{\theta}}}T\int_0^{1/2}
\min\left(1,\frac{1}{H\Vert u \Vert}\right)du\nonumber\\
&&\ll T^{3/2}H^{-1}\log H\ll T^{-1/2}{\cal L}.\nonumber
 \end{eqnarray}

So we have(noting trivially $G(x,H)\ll T^{1/2}$)
\begin{eqnarray}
\int_T^{2T}R_2^2(x,H)dx&&\ll T^{2l}+T^{2l-1}\int_T^{2T}G^2(x,H)dx\\
&&\ll T^{2l}+T^{2l-1/2}\int_T^{2T}G(x,H)dx\nonumber\\
&&\ll T^{2l}+T^{2l-1} {\cal L}\ll T^{2l}. \nonumber
\end{eqnarray}

\subsection{\bf \bf Proof of Theorem 1}

We have
\begin{eqnarray}
R^2(2\pi x)=R_1^2(x,H)+2R_1(x,H)R_2(x,H)+R_2^2(x,H).
\end{eqnarray}

From (5.26), (5.28) and Cauchy's inequality we have
\begin{eqnarray}
\int_T^{2T}R_1(x,H)R_2(x,H)dx\ll T^{2l+1/4}.
\end{eqnarray}
From (5.26), (5.28), (5.29) and (5.30) we get
\begin{eqnarray}
\int_T^{2T}R^2(2\pi
x)dx=\frac{2^{3-2l}C_{l,\theta}}{\pi^2(l-1)!^2}\int_T^{2T}x^{2l-\frac{1}{2}}dx
 +O(T^{2l +\frac{4\gamma+1}{8\gamma+4} +\varepsilon}  ).
\end{eqnarray}

Hence
\begin{eqnarray}
\int_1^{T}R^2(2\pi
x)dx&&=\frac{2^{3-2l}C_{l,\theta}}{\pi^2(l-1)!^2}\int_1^{T}x^{2l-\frac{1}{2}}dx
 +O(T^{2l +\frac{4\gamma+1}{8\gamma+4} +\varepsilon}     )\\
 &&=\frac{2^{3-2l}C_{l,\theta}}{(2l+1/2)\pi^2(l-1)!^2} T^{2l+1/2}
 +O(T^{2l +\frac{4\gamma+1}{8\gamma+4} +\varepsilon}   ). \nonumber
\end{eqnarray}
Now Theorem 1 follows from (5.32).

School of Mathematical Sciences,\\
Shandong Normal University,\\
Jinan, Shandong, 250014,\\
P.R.China\\
E-mail:zhaiwg@hotmail.com

\bigskip

\noindent

Current address:

Department of Mathematics,

China University of Mining and Technology(Beijing),

Beijing 100083, P. R. China


\begin{thebibliography}{99}

\bibitem{BG}V. Bentkus, F. G\"otze, Lattice point problems and
distribution of values of quadratic forms, Ann. of Math. (2){\bf
50}: 3(1999), 977-1027.

\bibitem{BE}P. H. B\'erard, On the wave equation on a compact
Riemannian  manifold without conjugate points, Math. Z. {\bf
155}:3(1977), 249-276.

\bibitem{BL}L. Bleher, On the distribution of the number of lattice
points inside a family of convex ovals, Duke Math. J. {\bf 67}:
3(1992), 461-481.



\bibitem{CPT}D. Chung, Y. N. Petridis and J. Toth, The
remainder in Weyl's law for Heisenberg manifolds II,  Bonner
Mathematische Schriften, Nr. 360, Bonn, 2003, 16 pages.



\bibitem{Fo}G. B. Folland, Harmonic Analysis in Phase Space,
Princeton University Press(1989), 9-73.

\bibitem{FR}F. Fricker, Einf\"uhrung in die Gitterpunketlehre, [Introduction to lattice point theory]
Lehrb\"ucher und Monographien aus dem Gebiete der Exakten
Wissenschaften(LMW), Mathematische Reihe[Textbooks and Monographs in
the Exact Sciences]73, Birkh\"auser Verlag, Basel-Boston, Mas.,
1982.

\bibitem{GW}C. Gordon , E. Wilson, The spectrum of the Laplacian on
Riemannian Heisenberg manifolds, Michigan Math. J. {\bf 33}(2)
(1986), 253-271.

\bibitem{G}F. G\"otze, Lattice point problems and
  values of quadratic forms, Inventiones Mathematicae{\bf
  157}(2004), 195-226.

\bibitem{Ha}G. H. Hardy, On the expression of a number as the sum of
two squares, Quart. J Math. {\bf 46}(1915), 263-283.



\bibitem{He2}D. R. Heath-Brown, The Piatetski-Shapiro prime theorem, J.of Number
theory, Vol{\bf 16}(1983), 242-266.

\bibitem{Ho}L. H\"ormander, The spectral function of an elliptic
operator, Acta Math.{\bf 121}(1968), 193-218.

\bibitem{Hu}M. N. Huxley, Exponential sums and lattice points III,
Proc. London Math. Soc. {\bf87} (3) (2003), 591--609.




\bibitem{IVr}V.YA. Ivrii, Precise Spectral  Asymptotics for elliptic
Operators Acting in Fibrings over Manifolds with Boundary, Springer
Lecture notes in Mathematics {\bf 1100}(1984).

\bibitem{Ka}I. K\'{a}tai, The number of lattice points in a circle (in
Russian), Ann. Univ. Sci. Budapest Rolando E\"{o}tv\"{o}s, Sect.
Math., 8 (1965), 39-60.

\bibitem{K}M. Khosravi, Third moment of the remainder error term in
Weyl's law for Heisenberg manifolds, arXiv: 0711.0073..

\bibitem{KP}M. Khosravi, Y. Petridis, The remainder in Weyl's law
for $n$-dimensional Heisenberg manifolds, Proc. of the American
Math. Soc. {\bf 133}(2005), 3561-3571.

\bibitem{KT}M. Khosravi, J. Toth, Cram\'er's formula for Heisenberg
manifolds, Ann. de l'institut Fourier {\bf 55}(2005), 2489-2520.

\bibitem{KH}A. Y. Khinchin,  Zur metrischen Theorie der diophantischen
Approximationen, Math. Z. {\bf 24} (1926), no. 4, 706-714.

\bibitem{Ku}M. Kuhleitner, W. G. Nowak, The asymptotic behaviour of
the mean-square of fractional part sums, Proc. Edinb. Math. Soc.
{\bf 43}(2000), 309-323.

\bibitem{Mi}S. H. Min, The methods of number theory(in Chinese),
Science Press, Beijing: 1981.



\bibitem{PT}Y. Petribis , J. Toth, The remainder in Weyl's law for
Heisenberg manifolds,  J. Differential Geom. {\bf 60}(2002),
455-483.


\bibitem{R}K. F. Roth, Rational aproximation to algebraic numbers,
Mathematica {\bf 2}(1955), 1-20.

\bibitem{St}E. M. Stein, Harmonic Analysis, Princeston
University Press(1993), 527-574.

\bibitem{Ts}Kai-Man Tsang, Higher-power moments of $\Delta(x), E(t)$ and $P(x)$,
Proc. London Math. Soc.(3){\bf 65}(1992), 65-84.

\bibitem{Va}J. D. Vaaler, Some extremal functions in Fourier analysis, Bull.
Amer. Math. Soc. {\bf 12} (1985), 183--216.

\bibitem{V} I. M. Vinogradov, Special variants of the method of trigonometric sums,
 (Nauka, Moscow), 1976; English transl. in his Selected works (Springer-Verlag),
 1985.

\bibitem{Vo}A. V. Volovoy, Improved two-term asymptotics for the
eigenvalue distribution  function of an elliptic operator on a
compact manifold, Comm. Partial Differemtial Equations{\bf 15}:
11(1990), 1509-1563.

\bibitem{Z1} Wenguang Zhai, On higher-power moments of $\Delta(x)$ (II), Acta Arith.
{\bf 114} (2004), 35--54.

\bibitem{Z2} Wenguang Zhai, On higher-power moments of $\Delta(x)$ (III), Acta Arith.
{\bf 118} (2005), 263--281.

\bibitem{Z4}Wenguang Zhai, On the error term in Weyl's law for the Heisenberg
manifolds, Acta Arith., in press(see also arXiv: 0805.3856).

\end{thebibliography}
\end{document}